\documentclass[11pt, a4paper]{article}
\usepackage{general}

\usepackage{fancyhdr}
\usepackage[
backend=biber,
style=alphabetic,
]{biblatex}
\DeclareFieldFormat{doi}{#1} 
\AtBeginDocument{\DeclareFieldFormat{doi}{}} 
\DeclareFieldFormat{url}{#1}
\AtBeginDocument{\DeclareFieldFormat{url}{}}
\DeclareFieldFormat{issn}{#1}
\AtBeginDocument{\DeclareFieldFormat{issn}{}}
\usepackage{subfiles}

\begin{document}

\thispagestyle{plain}

\vskip 3.0em     

\begin{center}
	{\bf \Large The 2-categorical S-matrix of a braided fusion 1-category
        \\[5pt]
        is a character table
	}

\vskip 2.6em

{\large
Alea Hofstetter\,$^{a},\,$
Christoph Schweigert\,$^{b}\,$
}

\vskip 2.6em

 \it$^a$
 Fachbereich Mathematik, \ Universit\"at Hamburg\\
 Bereich Algebra und Zahlentheorie\\
 Bundesstra\ss e 55, \ D\,--\,20\,146\, Hamburg
 \\[9pt]
 \it$^b$
 Fachbereich Mathematik, \ Universit\"at Hamburg\\
 Bereich Algebra und Zahlentheorie\\
 Bundesstra\ss e 55, \ D\,--\,20\,146\, Hamburg

\end{center}

\vskip 4.7em

\noindent{\sc Abstract}\\[3pt]
The semisimple module categories over a braided fusion category $\C$ form a connected fusion 2-category $\Mod(\C)$. Its Drinfeld center $\mathcal{Z}(\Mod(\C))$ is a 
braided fusion 2-category. To any braided fusion 2-category,
Johnson-Freyd and Reutter \cite{JFR23} have associated a
matrix-valued invariant, the 2-categorical $\St$-matrix.
In this short note we investigate this matrix of $\mathcal{Z}(\Mod(\C))$ as an invariant for the
braided fusion 1-category $\C$ and show that it reduces to 
the character table of the M\"uger center of $\C$.
  

\section{Introduction}

Braided monoidal categories are ubiquitous in quantum topology and
representation theory and have been intensely studied, particularly
in the case of fusion categories. Indeed, the classification of fusion
categories has been an active field of research for a long time,
making it essential to have computable invariants at hand to distinguish them.
In the case of modular fusion categories, modular data have been used for
decades, but the results of \cite{MS21} show that inequivalent modular
categories can have the same modular data.

Hence, it is reasonable to look for more invariants, for example as
follows: any braided fusion category $\C$ gives rise, via its semisimple module
categories, to a connected fusion 2-category and thus a four-dimensional
topological field theory. Its Drinfeld center is a braided fusion 2-category to
which Johnson-Freyd and Reutter \cite{JFR23} have associated a
matrix-valued invariant which we call the 2-categorical $\St$-matrix.
The question arises whether this is another interesting invariant for
the braided fusion category $\C$.

The answer we give in this short note for braided pivotal fusion categories
is, unfortunately, negative: 
we show in our main theorem \ref{theorem:main}
that the 2-categorical $\St$-matrix equals
the character table of the M\"uger center of $\C$. We confirm this
result by computing the $\St$-matrix for pointed braided fusion categories.
Along this way, we also 
gain some additional insight into the structure of braided
module categories and the structure of module braidings of the regular module
category of an arbitrary braided monoidal category.

We conclude this brief introduction by putting the result into context,
using some recent results in the Morita theory of fusion 2-categories. 
Morita theory is relevant here, since by \cite[Theorem 2.3.2]{Dec24}
Morita equivalent fusion 2-categories have the same Drinfeld centers.
It has been shown \cite[Theorem 4.1.6]{Dec24} that a fusion 2-category is Morita 
equivalent to the 2-Deligne product of a strongly fusion 2-category and an invertible fusion 2-category, i.e. a fusion 2-category of the form 
$\Mod(\C)$ with $\C$ being a non-degenerate braided fusion 1-category.
The latter does not contribute to the Drinfeld center 
\cite[Corollary 5.3.3]{Dec24} and can, thus, be disregarded for our purposes.
By \cite{JFY21} these categories are governed by a group
and a certain cocycle in a (super)cohomology theory. Unfortunately, 
this is not suitable as a direct starting point for computations: 
at present, it is difficult to make the cocycles explicit in the 
so-called fermionic case. For this reason, we decided to address the
problem in this paper by direct methods which are essentially 1-categorical.
Our results are in agreement with the much more comprehensive theory
presented in \cite{Dec24}.

In Section \ref{section2}, we review the relevant notions and propositions
and classify in Proposition \ref{braidingregular} module braidings on
the regular module category.
Section \ref{section31} presents the new results: the identification of
the 2-categorical $\St$-matrix with the character table of the
M\"uger center, cf.\ Theorem
\ref{theorem:main}. In Section 
\ref{section32}, we compute the $\St$-matrix for a
braided pointed fusion category explicitly.

\subsubsection*{Acknowledgments}

We are grateful to Theo Johnson-Freyd,
Hannes Kn\"otzele and David Reutter for discussions.
The authors acknowledge support by the Deutsche
Forschungsgemeinschaft (DFG, German Research Foundation) under Germany's
Excellence Strategy - EXC 2121 ``Quantum Universe'' - 390833306 and the
Collaborative Research Center - SFB 1624 ``Higher structures, moduli
spaces and integrability'' - 506632645.

\section{Preliminaries and definitions}\label{section2}

We work over $\K$, an algebraically closed field of characteristic zero. All categories are assumed to be $\K$-linear. 1-categories are denoted by calligraphic letters $\B,\C,\dots$ and 2-categories are denoted by fractured letters $\FB,\FC,\dots$.

A fusion 1-category is a rigid semisimple linear monoidal category with only finitely many isomorphism classes of simple objects such that the unit object is simple. A braiding on a fusion 1-category is a coherent natural isomorphism from the tensor product to the opposite tensor product. The notion of fusion 2-categories was introduced in \cite{DR18}. For a fusion 2-category we have the notion of a braiding as well. Though, the natural isomorphism is replaced by an equivalence and thus, we have some coherence data. For the general notion of a braided 2-category, we refer to \cite{Cr98}. We assume the reader to be familiar with the concept of (braided) fusion 1-categories as well as (braided) fusion 2-categories. The notions central to this paper will though be recalled. 

 Let $(\C,\I,\otimes, c)$ be a braided fusion 1-category. Here, $\I\in\C$ denotes the unit object, $\otimes$ the monoidal product and $(c_{X,Y}:X\otimes Y\to Y\otimes X)_{X,Y\in\C}$ denotes the braiding. Usually, we omit the unit object, the monoidal product and the braiding from our notation. Furthermore, we usually assume $\C$ to be strict. We will use graphical calculus as well to denote the braiding as follows. 
    \begin{center}
        \begin{tikzpicture}[scale=0.4, every node/.style={scale=0.6}]
            \draw[out=90, in=-90] (2,0) to (0,4);
            \fill[white] (1,2) circle (0.1);
            \draw[out=90, in=-90] (0,0) to (2,4);
            \node at (1, -1) {$c_{X,Y}$};
            \node at (0, -0.3) {$X$};
            \node at (2, -0.3) {$Y$};
        \end{tikzpicture}
        \hspace {3cm}
        \begin{tikzpicture}[scale=0.4, every node/.style={scale=0.6}]
            \draw[out=90, in=-90] (0,0) to (2,4);
            \fill[white] (1,2) circle (0.1);
            \draw[out=90, in=-90] (2,0) to (0,4);
            \node at (1, -1) {$c_{Y,X}^{-1}$};
            \node at (0, -0.3) {$X$};
            \node at (2, -0.3) {$Y$};
        \end{tikzpicture}
    \end{center}
The M\"uger center, also called the symmetric center, of a braided fusion category $\C$ is the full monoidal subcategory $\mathcal{Z}_2(\C)$ consisting of objects $X\in \C$ such that for all $Y\in\C$ the double braiding $c_{Y,X}\circ c_{X,Y}$ is trivial. The Drinfeld center of a (not necessarily braided) fusion 1-category $\C$ is a braided fusion 1-category with objects being pairs consisting of objects of $\C$ and half-braidings on those objects. We will denote it by $\mathcal{Z}(\C)$.

A braided monoidal 2-category is denoted by $(\FB,\I, \Box, \br)$ with $\I$ denoting the unit object, $\Box$ the tensor product and $\br$ the braiding. Vertical composition of 2-morphisms is denoted by $\cdot_v$ and horizontal composition by $\cdot_h$. Furthermore, we assume the underlying 2-category to be strict. Further assumptions concerning the strictness of the monoidal structure and the braiding will be specified where needed. The braiding consists, for every pair of objects $A,B\in \FB$, of 1-morphisms $\br_{A,B}: A\Box B\to B\Box A$ which assemble to a natural equivalence compatible with the tensor product up to modifications. Additionally, we have for all 1-morphisms $f:A\to A'$ and $g:B\to B'$ 2-isomorphisms filling the naturality squares: \begin{equation*} \begin{tikzpicture}[scale=1.2, every node/.style={scale=0.8}]
    \node at (0,2) {$A\Box B$};
    \node at (2,2) {$B\Box A$};
    \node at (0,0) {$A'\Box B$};
    \node at (2,0) {$B\Box A'$};
    \draw[->] (0.5,0) to (1.5,0);
    \node at (1,-0.3) {$\br_{A',B}$};
    \draw[->] (0.5,2) to (1.5,2);
    \node at (1,2.3) {$\br_{A,B}$};
    \draw[->] (0,1.5) to (0,0.5);
    \node at (-0.4,1) {$f\Box B$};
    \draw[->] (2,1.5) to (2,0.5);
    \node at (2.4,1) {$B\Box f$};
    \draw[->, double] (0.8,0.8) to (1.2,1.2);
    \node at (0.7, 1.3) {$\br_{f,B}$};
\end{tikzpicture}
\hspace{2cm}
\begin{tikzpicture}[scale=1.2, every node/.style={scale=0.8}]
    \node at (0,2) {$A\Box B$};
    \node at (2,2) {$B\Box A$};
    \node at (0,0) {$A\Box B'$};
    \node at (2,0) {$B'\Box A$};
    \draw[->] (0.5,0) to (1.5,0);
    \node at (1,-0.3) {$\br_{A,B'}$};
    \draw[->] (0.5,2) to (1.5,2);
    \node at (1,2.3) {$\br_{A,B}$};
    \draw[->] (0,1.5) to (0,0.5);
    \node at (-0.4,1) {$A\Box g$};
    \draw[->] (2,1.5) to (2,0.5);
    \node at (2.4,1) {$g\Box B$};
    \draw[->, double] (0.8,0.8) to (1.2,1.2);
    \node at (0.7, 1.3) {$\br_{A,g}$};
\end{tikzpicture}
\end{equation*}

For a monoidal 2-category $\FC$ there is the notion of a Drinfeld center $\mathcal{Z}(\FC)$ defined similarly as in the setting of 1-categories (compare \cite[Section 3.1.1]{Cr98}, \cite[Section 2.3]{DN21}).

Two simple objects of a fusion 1-category are, by Schur's lemma, either isomorphic or admit only a zero morphism. This implies as well that for all simple objects $X\in\C$ every endomorphism $f:\;X\to X$ is defined by a scalar in $\K$. We denote this scalar by $\langle f\rangle\in\K$. For a fusion  2-category, the situation is more involved:\; non-equivalent objects can admit non-trivial morphisms. This leads to an equivalence relation that is coarser:\;
\begin{center}
    $X\sim Y$ if there is a non-zero 1-morphism $f:\;X\to Y$
\end{center}
(compare \cite[Theorem 2.7]{JFR23}). The corresponding equivalence classes are called components or Schur equivalence classes of $\FC$. We denote the set of components by $\pi_0(\FC)$. 
In this paper, we will study a matrix assigned to any braided fusion 2-category $\FB$, called the 2-categorical $\St$-matrix; the elements of the set $\pi_0(\FB)$ label its rows. 

A fusion 2-category with one component is called connected. It can be obtained from a braided fusion 1-category $\C$ as follows. Take the fusion 2-category $\Mod(\C)$ of the (right) finite semisimple module categories, or equivalently the separable algebra objects \cite[Construction 2.1.19]{DR18}. Conversely, let $\FC$ be a fusion 2-category with monoidal unit $\I$, then $\End_\FC(\I)$ is a braided fusion 1-category. In fact, all connected fusion 2-categories arise as the category of modules of some braided fusion 1-category:

\begin{proposition}\cite[Remark 2.1.22]{DR18}
	Let $\FC$ be a connected fusion 2-category. Then the 1-category of endomorphisms of the unit object $\End_\FC(\I)$ is a braided fusion 1-category and $$\FC\overset{\otimes}{\simeq}\Mod(\End_\FC(\I)).$$    
\end{proposition}

We will focus in the sequel on a specific class 
of braided fusion 2-categories: the 2-categorical Drinfeld 
$\mathcal{Z}(\FC)$ of a connected fusion 2-category $\FC$,
which is a (braided) fusion 2-category by \cite[Lemma 2.18]{JFR23}.

 In the case of a connected fusion 2-category $\FC=\Mod(\C)$, for $\C$ a braided fusion 1-category, the 2-categorical Drinfeld center can be described in two equivalent ways:\; via separable half-braided algebras and via finite semisimple braided module categories. A separable half-braided algebra $(A,\gamma)$ is a separable algebra $A$ in $\C$ together with the additional structure of an additional half-braiding $\gamma$ obeying consistency conditions with the multiplication and braiding that can be found in \cite[Section 2.4]{JFR23}. Together with bimodules obeying a compatibility condition with the braiding and the additional half-braiding, they form a braided fusion 2-category $\sHBA(\C)$. For more details, we refer to \cite[Section 2.4]{JFR23}. An equivalent braided fusion 2-category is the 2-category $\BrMod(\C)$ of finite semisimple braided (right) $\C$-module categories whose objects are finite semisimple $\C$-module categories $\M$ together with a natural isomorphism $(\sigma_{M,X}:\;M\ta X\to M\ta X)_{X\in \C,M\in \M}$ (compare \cite[Section 2.2]{JFR23}, \cite[Section 4]{DN21}). The module braiding $\sigma$ has to satisfy that \begin{equation}\label{braidedmodulecat0}
     \sigma_{M, \I}=\id_M
 \end{equation} for all $M\in \M$ and the following two compatibility conditions with the braiding and the tensor product on $\C$: \begin{equation}\label{braidedmodulecat1}
     \begin{tikzcd}
            (M\ta X) \ta Y\arrow[rr, "\sigma_{M\ta X,Y}"] \arrow[d, "m^{-1}_{M,X,Y}"'] && (M\ta X) \ta Y \arrow[d, "m^{-1}_{M,X,Y}"]\\
             M \ta (X\otimes Y) \arrow[d, "\id_M\ta c_{X,Y}"'] && M \ta (X\otimes Y) \arrow[d, "\id_M\ta c_{Y,X}^{-1}"]\\
            M \ta (Y\otimes X) \arrow[d, "m_{M,Y,X}"'] && M \ta (Y\otimes X) \arrow[d, "m_{M,Y,X}"]\\
            (M\ta Y)\ta X \arrow[rr, "\sigma_{M,Y}\ta \id_X"]&& (M\ta Y)\ta X
        \end{tikzcd}
\end{equation}
\begin{equation}\label{braidedmodulecat2}
    \begin{tikzcd}
             M\ta (X\otimes Y)\arrow[rr, "\sigma_{M, X\otimes Y}"] \arrow[d, "\id_M \ta c_{Y,X}^{-1}"'] && M\ta (X\otimes Y) \arrow[d, "\id_M\ta c_{X,Y}"]\\
            M\ta (Y\otimes X) \arrow[d, "m_{M,Y, X}"'] && M\ta (Y\otimes X) \arrow[d, "m_{M,Y,X}"]\\
            (M\ta Y) \ta X \arrow[rd, "\sigma_{M\ta Y, X}"']&& (M\ta Y) \ta X\\
            & (M\ta Y) \ta X\arrow[ru, "\sigma_{M, Y}\ta \id_X"']
        \end{tikzcd}
\end{equation}

The regular module category $\C_\C$ of a braided fusion 1-category
can always be equipped with a module braiding defined by the 
double braiding $c$ on $\C$, namely $\sigma_{Y, X}=c_{X,Y}\circ c_{Y,X}:\; Y\ta X\to Y\ta X$. 
For later use, we describe all module braidings on the regular module category.

\begin{proposition}\label{braidingregular}
Let $\C$ be a braided monoidal category. (Here we need not to
make any additional assumption like linearity or finiteness.) Then
any monoidal automorphism $(\eta_X)_{X\in\C}$ on the identity
functor on $\C$, with its trivial structure of a monoidal functor, provides
a module braiding on the regular $\C$-module category via
\begin{equation}\label{braidingregulareq}
    \sigma_{Y, X}\coloneqq c_{X,Y}\circ c_{Y,X}\circ (Y\otimes \eta_{X}): Y\otimes X\to Y\otimes X.
\end{equation}
Conversely, any module braiding on the regular module category is
of this form.
\end{proposition}

\begin{proof}
We can assume that $\C$ is strict.
For any monoidal automorphism $\eta$ of the identity functor, we have 
$\eta_\I=\id_\I$, hence
equation (\ref{braidedmodulecat0}) is automatically fulfilled.
A straightforward computation shows that (\ref{braidingregulareq}) indeed
defines a module braiding fulfilling \eqref{braidedmodulecat1} and \eqref{braidedmodulecat2}.

Conversely, suppose $\sigma$ is a module braiding on the 
regular module category $\C_\C$. Then combining 
the conditions \eqref{braidedmodulecat1} and \eqref{braidedmodulecat2}
yields
\begin{equation}\label{braidedmoduleregular2}
        \sigma_{\I, X\otimes Y}=  (X\otimes \sigma_{\I,Y}) \circ (\sigma_{\I,X}\otimes Y) =(\sigma_{\I, X}\otimes Y)\circ (X\otimes \sigma_{\I,Y})
    \end{equation}
for all objects $X,Y\in \C$.   
The isomorphisms $\sigma_{\I,X}: X\cong \I\otimes X \to \I\otimes X\cong X$   
are natural and form a monoidal automorphism due to 
(\ref{braidedmoduleregular2}). 
\end{proof}

\begin{example}\label{regularmodulebraidings}
Consider $\C=\vect_G$, the abelian category of finite-dimensional 
$G$-graded $\Cn$-vector spaces for some finite abelian group $G$. Braidings 
and associators on $\vect_G$ with the underlying tensor product functor being 
the $G$-graded tensor product of vector spaces can be classified via abelian 3-cocycles $[(\Psi, \Omega)]\in H^3_{ab}(G,\Cn^*)$ or equivalently via quadratic forms $q:\;G\to \K^\times$ (compare  \cite[Appendix E]{MS89} and \cite{TFT3}). We fix an abelian three-cocycle $(\Psi, \Omega)$ defining the braiding on $\vect_G$.

By Proposition \ref{braidingregular}, we need to classify monoidal 
automorphisms $\eta$ of the identity monoidal functor. They are
given by automorphisms $\eta_g:\so{g}\to\so{g}$ which are scalar multiples
of the identity $\eta_g=\varphi(g)\id_{\so{g}}$. By monoidality,
$\varphi$ is a character on $G$.

Later on we will show that these are, up to Schur equivalence, all finite semisimple indecomposable braided $\vect_G$-module categories.
\end{example}

\begin{proposition}\label{BrModHBAequ}
    For a braided fusion 1-category $\C$ we have $$\mathcal{Z}(\Mod(\C))\simeq \BrMod(\C)\simeq\sHBA(\C)$$ as braided fusion 2-categories.    
\end{proposition}

For a proof of the first part we refer to \cite[Theorem 4.11]{DN21} and for the second part we refer to \cite[Theorem 2.40]{JFR23}.
    
As already mentioned above an $\St$-matrix has been defined for a braided fusion 2-category, compare \cite[Subsection 2.7]{JFR23}. We recall the definition and the most relevant computational concepts for the convenience for the reader.

\begin{construction}   
Let $\FB$ be a braided fusion 2-category. Let $A\in \FB$ be a simple object and let $b:\;\I\to\I$ be a simple 1-endomorphism of the unit object. The braiding on $\FB$ gives 2-isomorphisms 
	\begin{eqnarray*}
		\br_{A,b}:\;\br_{A,\I}\circ(A\Box b)\Rightarrow (b\Box A)\circ \br_{A,\I} \,\,\, ,\\
		\br_{b,A}:\;\br_{\I,A}\circ(b\Box A)\Rightarrow (A \Box  b)\circ \br_{\I,A} \,\,\,.      
	\end{eqnarray*}
With the strictness we impose, we have $\I\Box A=A=A\Box\I$ and that
 $\br_{\I,A}$ and $\br_{A,\I}$ are trivial. We thus obtain 2-isomorphisms \begin{eqnarray*}
		\br_{A,b}:\;A\Box b\Rightarrow b\Box A,\\
		\br_{b,A}:\;b\Box A \Rightarrow A\Box b.
	\end{eqnarray*}
	Now, we can define the full braid as the composite
    $$\br_{b,A}\cdot_v\br_{A,b}:\;A\Box b\Rightarrow A\Box b.$$
	
	To obtain  scalars, assume that $b\in \Hom_\FB(\I,\I)$ is simple. Furthermore, choose duality data 
    \begin{eqnarray*}
		\ev_b:\;b\circ b^*\to\id_\I,\\
		\coev_b:\;\id_\I\to b^*\circ b
	\end{eqnarray*}
    and analogously other sided duals $^*b$ as well as a 2-isomorphism $\alpha_b:\; {^{*}{b}}\Rightarrow b^*$. Note that $\alpha_b$ only necessarily exists for $\Hom_\FB(\I,\I)$ being semisimple \cite[Proposition 4.8.1]{EGNO}. Let $d_{\alpha_b}$ be the scalar defining the 2-endomorphism $\ev_b\cdot_v(\id_{b} \Box \alpha_b)\cdot_v\coev_{{^{*}{b}}}$. Now, we define $R_{A,b}:\;\id_A\Rightarrow\id_A$ to be the 2-endomorphism \begin{equation}\label{RTilde}
		\frac{1}{d_{\alpha_b}}((A\Box\ev_b)\cdot_v((\br_{b,A}\cdot_v\br_{A,b})\Box\alpha_b)\cdot_v(A\Box\coev_{{^{*}{b}}})).
	\end{equation}
    Since $A\in\mathfrak{B}$ is simple, the 1-morphism $\id_A\in \End_\FB(A)$
    is simple and thus, $R_{A,b}$ is a multiple of the identity as a 2-morphism by Schur's lemma. We denote this scalar by $\St_{A,b}$.
\end{construction}

\begin{lemma}
	The scalar $\St_{A,b}$ depends only on the Schur equivalence
    class of $A\in\mathfrak{B}$. Furthermore, it is independent of the representative of the isomorphism class $b\in \Hom_\FB(\I,\I)$.  
\end{lemma}

\begin{proof}
	Assume $f:\;A\to B$ is a non-zero 1-morphism in $\mathfrak{B}$. The naturality of the braiding (see \cite[p.194]{Cr98}) together with $\br_{\I,\_}$, $\br_{\_,\I}$ being trivial and the interchange 2-isomorphism of the monoidal structure of $\FB$ respecting identities imply that $$\br_{b,A}=\br_{b,B}.$$ Furthermore, let $b':\;\I\to\I$ be another simple object that lies in the same equivalence class, i.e. there is a non-zero 2-isomorphism $\alpha:\;b\Rightarrow b'$. Then $$(A\Box \alpha)\cdot_v\br_{b,A}=\br_{b',A}\cdot_v(\alpha\Box A)$$ by naturality (see \cite[Lemma 7, 2nd condition]{BN96}) as well.
\end{proof}

\begin{lemma}
	The scalar $\St_{A,b}$ is independent of the choice of 2-isomorphism $\alpha_b:\;{^{*}{b}}\Rightarrow b^*$.
\end{lemma}

\begin{proof}
    This follows analogously to the fact that the trace is well-defined for fusion 1-categories. We spell it out in the following. Let $\alpha_b:\;{^{*}{b}}\Rightarrow b^*$ be some choice of 2-isomorphism. Assume that $\alpha'_b:\;{^{*}{b}}\Rightarrow b^*$ is another choice of 2-isomorphism. Then $\alpha'_b=\alpha_b\cdot_v\alpha_b^{-1}\cdot_v\alpha'_b$, and $\alpha_b^{-1}\cdot_v\alpha'_b:\;{^{*}{b}}\to {^{*}{b}}$ is a 2-automorphism. Thus, for fixed $\alpha_b:\;{^{*}{b}}\Rightarrow b^*$, any other choice of 2-isomorphism $\alpha_b'$ is defined by some 2-isomorphism of a simple object $\psi_b:\;{^{*}{b}}\to {^{*}{b}}$ via $\alpha_b'=\alpha_b\cdot_v \psi_b$. Furthermore, $\psi_b$ is a morphism of simple objects in $\End_\FB({\I})$ and thus defined by some scalar. Therefore, $$d_{\alpha_b'}=d_{\alpha_b\cdot_v \psi_b}=\langle\ev_b\cdot_v(\id_{b} \Box (\alpha_b\cdot_v \psi_b))\cdot_v\coev_{{^{*}{b}}}\rangle=\langle\psi_b\rangle \cdot d_{\alpha_b}.$$ Moreover, we have \begin{align*}&\langle(A\Box\ev_b)\cdot_v((\br_{b,A}\cdot_v\br_{A,b})\Box(\alpha_b\cdot_v\psi_b))\cdot_v(A\Box\coev_{{^{*}{b}}})\rangle\\    &\hspace{1cm}=\langle\psi_b\rangle\cdot\langle(A\Box\ev_b)\cdot_v((\br_{b,A}\cdot_v\br_{A,b})\Box\alpha_b)\cdot_v(A\Box\coev_{{^{*}{b}}})\rangle.
	\end{align*} Thus, $\St_{A,b}$ is independent of the choice of 2-isomorphism $\alpha_b$ as $\langle\psi_b\rangle$ occurs in the nominator as well as in the denominator of the definition of $\St_{A,b}$. 
\end{proof}

In regard of the last two lemmas, we are ready to define the 
2-categorical $\St$-matrix.

\begin{definition} (\cite[Lemma 2.56]{JFR23})
Let $\FB$ be a braided fusion 2-category.
	The 2-categorical $\St$-matrix is the function
	$$\St:\;\pi_0\FB\times\pi_0\Hom_\FB(\I,\I)\to\K:\; (A,b)\mapsto \St_{A,b}.$$
	Here, $\pi_0\mathfrak{B}$ denotes the set of Schur equivalence classes of the 2-category $\mathfrak{B}$ while $\pi_0\Hom_\FB(\I,\I)$ denotes the set of isomorphism classes of objects in the fusion 1-category $\Hom_\FB(\I,\I)$. 
\end{definition}

We now turn to the aforementioned special class of braided fusion 2-categories:\; to the case of a Drinfeld center of a connected fusion 2-category, i.e. $\FB=\mathcal{Z}(\Mod(\C))$ for some braided fusion 1-category $(\C,c)$. In this case, we can compute the $\St$-matrix pairing more explicitly. Recall the following proposition about the relation of the endomorphism space of the unit object of $\mathcal{Z}(\Mod(\C))$ and the M\"uger center of $\C$. The proposition can be found in \cite[Lemma 2.16]{JFR23}.

\begin{proposition}
    For a braided fusion 1-category $\C$, $\Hom_{\mathcal{Z}(\Mod(\C))}(\I,\I)$ is equivalent to $\mathcal{Z}_2(\C)$ as braided fusion 1-categories.    
\end{proposition}

\begin{proof}
       In regard of Theorem \ref{BrModHBAequ}, rather than constructing the equivalence of $\Hom_{\mathcal{Z}(\Mod(\C))}(\I,\I)$ and $\mathcal{Z}_2(\C)$, we will construct an equivalence of $\Hom_{\BrMod(\C)}(\I,\I)$ and $\mathcal{Z}_2(\C)$. We consider $$({\C_\C}, (\sigma_{Y,X}=c_{X,Y}\circ c_{Y,X})_{X\in\C, Y\in {{\C}_{\C}}})$$ the monoidal unit in $\BrMod(\C)$. A braided module morphism $$F:({\C_\C}, \sigma)\to({\C_\C}, \sigma)$$ is in particular a module functor of the regular $\C$-module category $F:{\C_\C}\to{\C_\C}$ and thus defined by some object $A\in\C$, $F(\I_\C)=A$. For $F$ to be a braided module functor, we must have for all $X\in\C$ \begin{alignat*}{2}
          &&~F(\sigma_{\I_\C, X})&=\sigma_{F(\I_\C),X}\\
          \Leftrightarrow &&~F(c_{X,\I_\C}\circ c_{\I_\C,X})&=c_{X, F(\I_\C)}\circ c_{F(\I_\C),X}\\
          \Leftrightarrow&&~\id_{A\otimes X}&=c_{X,A}\circ c_{A,X}.
      \end{alignat*} Thus, $F$ is a braided module functor if and only if $A\in\mathcal{Z}_2(\C)$. This establishes the equivalence.
    
\end{proof}

Using the  equivalence in the previous proposition for one
index of the S-matrix and the two descriptions of $\FB=\mathcal{Z}(\Mod(\C))$ via separable half-braided algebras in $\C$ and finite semisimple braided module categories on $\C$ respectively for the other index, one can state two explicit formulae to compute the 2-categorical $\St$-matrix of $\FB$.

\begin{observation} \label{halfbraidedalgebrabraidedmoduleSmatrixcomputation}\cite[Example 2.54, Example 2.55]{JFR23}
    First, consider $(A,\gamma)\in\sHBA(\C)$ some separable half-braided algebra in $\C$, simple object in $\sHBA$, and simple $b\in \mathcal{Z}_2(\C)$ considered as a simple object in $\Hom_{\sHBA(\C)}(\I,\I)$. Choose some duality data for $b$ and some isomorphism $\alpha_b:\;{^*b}\to b^*$. Then we have \begin{align*}
        \St_{((A,\gamma),b)}
        &=\frac{1}{d_{\alpha_b}}\langle((\id_A\otimes\ev_b)\cdot((\gamma_{b,A}\circ c_{A,b})\otimes\alpha_b)\circ(\id_A\otimes\coev_{{^{*}{b}}})\rangle
    \end{align*} 
    Here, we can apply Schur's Lemma as we can see $(\id_A\otimes\ev_b)\cdot((\gamma_{b,A}\circ c_{A,b})\otimes\alpha_b)\circ(\id_A\otimes\coev_{{^{*}{b}}})$ as a 2-morphism $_AA_A\to {_AA}_A$ with $_AA_A$ seen as a half-braided bimodule and thus a 1-morphism $(A,\gamma)\to (A, \gamma)$ in $\sHBA$. Graphically, we can denote this 2-morphism as follows.
    \begin{center}
        \begin{tikzpicture}[scale=2.3, rotate=180, every node/.style={scale=0.8}]
        \node at (1.2, 0.5) {$\frac{1}{d_{\alpha_b}}$};
        \node at (0.5,1.1) {$A$};
        \node at (0.125, 0.9) {$b$};
        \draw[in=0, out=90] (0.75, 0.5) to (0.25, 0.75);
        \draw[in=90, out=180] (0.25, 0.75) to (0,0.5);
        \draw[in=180, out=-90] (0,0.5) to (0.25, 0.25);
        \draw[white, line width=1mm] (0.5,0.5) to (0.5,1);
        \draw (0.5,0) to (0.5,1);
        \draw[out=0, in=-90] (0.25, 0.25) to (0.75, 0.5);
        \node at (0.7, 0.2) {$\gamma_{b, A}$};
        \node at (0.7, 0.8) {$c_{A, b}$};
        \fill (0.475, 0.25) rectangle (0.525, 0.3);
        \fill (0,0.5) circle (0.02);
        \node at (-0.2, 0.5) {$\alpha_b$};
        \end{tikzpicture}  
    \end{center}
Note that we use the by $\gamma$ defined half braiding on the half-braided bimodule $_AA_A$ denoted by $\gamma_{\_,A}$.

For the other description of objects in $\FB$,    
consider $(\M,\sigma)\in \BrMod(\C)$, a finite semisimple braided module category over $\C$, simple object in $\BrMod(\C)$. As above, let $b\in \mathcal{Z}_2(\C)$ be simple considered as a simple object in $\Hom_{\sHBA(\C)}(\I,\I)$. Choose again some duality data for $b$ and some isomorphism $\alpha_b:\;{^*b}\to b^*$. Then we can compute the 2-categorical $\St$-matrix for any simple $M\in\M$ via \begin{equation*}
    \begin{split}
        \St_{((\M,\sigma),b)}=\frac{1}{d_{\alpha_b}}\langle(M\cong M\triangleleft I\overset{(\id_M\triangleleft\coev_b)}{\longrightarrow} M\triangleleft(b\otimes b^*) \cong (M\triangleleft b)\triangleleft b^* \overset{\sigma_{M,b}\triangleleft \id_b^*} {\longrightarrow} \\  (M\triangleleft b)\triangleleft b^* \cong M\triangleleft(b\otimes b^*) \overset{(\id_M\triangleleft\ev_b)}{\longrightarrow} M\triangleleft I\cong M)\rangle
    \end{split}
\end{equation*} The value of $\St_{((\M,\sigma),b)}$ is constant for all simple $M\in\M$.
\end{observation}

We will use the latter formula in concrete examples. For
the general case, we note that
instead of computing the 2-categorical $\St$-matrix explicitly via braided module categories or half-braided algebras, one can follow the proof of \cite[Theorem 2.57]{JFR23} and consider instead (parts of) the 1-categorical $\St$-matrix of the Drinfeld center of $\C$. 

\begin{observation}\label{Smatrixcomputation1categoryviaDrinfeldcenter}
Let $\mathcal{Z}(\C)$ denote the Drinfeld center of the braided fusion 1-category $\C$. The Drinfeld center is a braided fusion 1-category with $\Sone$-matrix denoted by $\Sone^{\mathcal{Z}(\C)}$. The 2-categorical Drinfeld center $\mathcal{Z}(\Mod(\C))$ gives a braided fusion 2-category. Its $\St$-matrix is denoted by $\St^{\mathcal{Z}(\Mod(\C))}$. The proof of \cite[Theorem 2.57]{JFR23} states that the 2-categorical 
$\St^{\mathcal{Z}(\Mod(\C))}$ is obtained from the 1-categorical $\Sone^{\mathcal{Z}(\C)}$ in two steps:\; \begin{enumerate}
        \item The rows and columns of $\Sone^{\mathcal{Z}(\C)}$ are labeled respectively by simple objects of $\mathcal{Z}(\C)$. We can embed $\mathcal{Z}_2(\C)$ into $\mathcal{Z}(\C)$ by $$X\mapsto (X,c_{\_,X})$$ for $X\in \mathcal{Z}_2(\C)$. Note that for the M\"uger center we only have one canonical embedding of $\mathcal{Z}_2(\C)\hookrightarrow \mathcal{Z}(\C)$ as the braiding $c_{\_,X}$ and the inverse braiding $c_{X,\_}^{-1}$ 
        coincide for $X\in \mathcal{Z}_2(\C)$. Thus, it makes sense to consider the $ |\pi_0(\mathcal{Z}(\C))|\times |\pi_0(\mathcal{Z}_2(\C))| $ matrix $\Sone^{\mathcal{Z}(\C)}_{\mathcal{Z}(\C), \mathcal{Z}_2(\C)}$, i.e. $\Sone^{\mathcal{Z}(\C)}$ where we only consider the columns indexed by simple objects in $\mathcal{Z}_2(\C)$. We obtain a matrix of rank $|\pi_0(\mathcal{Z}_2(\C))|$ as $\mathcal{Z}(\C)$ is a non-degenerate fusion 1-category.
        \item  The matrix $\Sone^{\mathcal{Z}(\C)}_{\mathcal{Z}(\C), \mathcal{Z}_2(\C)}$ has $|\pi_0(\mathcal{Z}_2(\C))|$ linearly independent rows. The other rows are just copies of these rows. By deleting the redundant rows, we obtain $\St^{\mathcal{Z}(\Mod(\C))}$. This can be seen as $\mathcal{Z}(\C)$ is a generator of $\BrMod(\C)$. For further details we refer to the proof of \cite[Theorem 2.57]{JFR23}. 
    \end{enumerate}
    Thus, instead of computing $\St^{\mathcal{Z}(\Mod(\C))}$ explicitly by considering Schur equivalence classes of braided module categories or half-braided algebras, we might as well consider $\Sone^{\mathcal{Z}(\C)}_{\mathcal{Z}(\C), \mathcal{Z}_2(\C)}$ and find $|\pi_0(\mathcal{Z}_2(\C))|$ many linearly independent rows. As linearly dependent rows are, in this special case, simply identical, we might as well just find $|\pi_0(\mathcal{Z}_2(\C))|$ many pairwise
    different rows. Those form the $\St^{\mathcal{Z}(\Mod(\C))}$-matrix.
\end{observation}

\section{Main result and examples}
In this section we compute the 2-categorical $\St$-matrix associated to a braided pivotal fusion 1-category $\C$. We will obtain the character table of the group associated to M\"uger center of $\C$. We confirm this result by computing the 2-categorical $\St$-matrix in a class of examples where we classify braided module categories up to Schur equivalence explicitly.

\subsection{Main result} \label{section31}

Let $\Rep(G)$ denote the category of finite-dimensional $\K$-representations of $G$ considered as a braided monoidal category with trivial braiding. Let $z\in G$ be a central element of order two, then $\Rep(G,z)$ denotes the same monoidal category, though the braiding is the one defined by $z$, compare \cite[Example 9.9.1(3)]{EGNO}. Recall that every symmetric fusion 1-category $\C$ is braided monoidal equivalent to either $\Rep(G)$ or $\Rep(G,z)$ for some finite group $G$ and $z\in G$ a central element of order two, see \cite[Theorem 9.9.22, Theorem 9.9.26]{EGNO}. In the first case $\C$ is called Tannakian and in the second case $\C$ is called super-Tannakian.

\begin{theorem}\label{theorem:main}
    Let $\C$ be a pivotal braided fusion 1-category with the M\"uger
    center being braided equivalent to $\Rep(G)$ or 
    $\Rep(G,z)$ with $G$ a finite group and $z\in G$ a central element of order two. Then the 2-categorical $\St$-matrix of $\mathcal{Z}(\Mod(\C))$ is, up to permutation of columns and rows, the character table of $G$.
\end{theorem}

\begin{proof}
    Following Observation \ref{Smatrixcomputation1categoryviaDrinfeldcenter} we  compute the 2-categorical $\St$-matrix of $\mathcal{Z}(\Mod(\C))$ via the 1-categorical $\St$-matrix of $\mathcal{Z}(\C)$. Assume $\C$ to be strict. Consider the inclusion $\tau:\;\mathcal{Z}_2(\C)\xhookrightarrow{} \mathcal{Z}(\C)$  from the M\"uger center into the Drinfeld center provided by the braiding on $\C$. Let $\mathcal{T}$ be the subcategory of $\mathcal{Z}(\C)$ defined by $\tau$. Let $z\in\mathcal{Z}(\C)$ and let $t\in\mathcal{T}$ be simple objects. The entry in the 1-categorical $\St$-matrix of $\mathcal{Z}(\C)$ in the column labeled by $t$ and the row labeled by $z$ are defined as follows. 
    \begin{center}
        \begin{tikzpicture}[scale=2.3, rotate=180, every node/.style={scale=0.8}]
        \node at (0.5,1.1) {$z$};
        \node at (0.125, 0.9) {$t$};
        \draw[in=0, out=90] (0.75, 0.5) to (0.25, 0.75);
        \draw[in=90, out=180] (0.25, 0.75) to (0,0.5);
        \draw[in=180, out=-90] (0,0.5) to (0.25, 0.25);
        \draw[white, line width=1mm] (0.5,0.5) to (0.5,1);
        \draw (0.5,0) to (0.5,1);
        \draw[out=0, in=-90, white, line width=1mm] (0.25, 0.25) to (0.75, 0.5);
        \draw[out=0, in=-90] (0.25, 0.25) to (0.75, 0.5);
        \node at (-0.1,0.5) {$\coloneqq$};
            \node at (-0.4,0.5) {$S_{z,t}$};
            \draw (-0.6,0) to (-0.6,1);
            \node at (-0.6,1.1) {$z$};
        \end{tikzpicture} 
    \end{center}
    As we assumed $\C$ to be pivotal, we omit the isomorphism $\alpha_t$ and the factor $d_{\alpha_t}$. Furthermore, the braiding depicted is the braiding in $\mathcal{Z}(\C)$. Let $t,t'\in\mathcal{T}$ be simple objects. Then we have by monoidality of the braiding in $\mathcal{Z}(\C)$ and by the pivotal structure the following.
    \begin{equation}\label{multiplicativity}
\begin{tikzpicture}[scale=2.3, rotate=180, every node/.style={scale=0.8}]
        \node at (0.5,1.1) {$z$};
        \node at (0.125, 0.9) {$t\otimes t'$};
        \draw[in=0, out=90] (0.75, 0.5) to (0.25, 0.75);
        \draw[in=90, out=180] (0.25, 0.75) to (0,0.5);
        \draw[in=0, out=90] (0.8, 0.5) to (0.25, 0.8);
        \draw[in=90, out=180] (0.25, 0.8) to (-0.1,0.5);
        \draw[in=180, out=-90] (0,0.5) to (0.25, 0.25);
        \draw[in=180, out=-90] (-0.1,0.5) to (0.25, 0.2);
        \draw[white, line width=1mm] (0.5,0.5) to (0.5,1);
        \draw (0.5,0) to (0.5,1);
        \draw[out=0, in=-90, white, line width=1mm] (0.25, 0.25) to (0.75, 0.5);
        \draw[out=0, in=-90] (0.25, 0.25) to (0.75, 0.5);
        \draw[out=0, in=-90, white, line width=1mm] (0.25, 0.2) to (0.8, 0.5);
        \draw[out=0, in=-90] (0.25, 0.2) to (0.8, 0.5);        
        \end{tikzpicture} 
        \begin{tikzpicture}[scale=2.3, rotate=180, every node/.style={scale=0.8}]
        \node at (1.2,0.5) {$=$};
        \node at (1,0.5) {$S_{z,t'}$};
        \node at (0.5,1.1) {$z$};
        \node at (0.125, 0.9) {$t$};
        \draw[in=0, out=90] (0.75, 0.5) to (0.25, 0.75);
        \draw[in=90, out=180] (0.25, 0.75) to (0,0.5);
        \draw[in=180, out=-90] (0,0.5) to (0.25, 0.25);
        \draw[white, line width=1mm] (0.5,0.5) to (0.5,1);
        \draw (0.5,0) to (0.5,1);
        \draw[out=0, in=-90, white, line width=1mm] (0.25, 0.25) to (0.75, 0.5);
        \draw[out=0, in=-90] (0.25, 0.25) to (0.75, 0.5); 
            \node at (-0.2,0.5) {$=$};
            \node at (-0.6,0.5) {$S_{z,t'}\cdot S_{z,t}$};
            \draw (-1,0) to (-1,1);
            \node at (-1,1.1) {$z$};            
        \end{tikzpicture}         
    \end{equation}

   If $\mathcal{Z}_2(\C)$ is Tannakian, let $G$ be a finite group such that $\Rep(G)\overset{\otimes, br}{\simeq}\mathcal{Z}_2(\C)$. Let $\tilde{\tau}:\;\Rep(G)\xhookrightarrow{} \mathcal{Z}(\C)$ be the braided monoidal inclusion given by the previous braided monoidal equivalence and the inclusion of the M\"uger center into the Drinfeld center. In the super-Tannakian case let $z\in G$ be a central element of order two such that $\Rep(G,z)\overset{\otimes, br.}{\simeq}\mathcal{Z}_2(\C)$ and let $\tilde{\tau}:\;\Rep(G,z)\xhookrightarrow{} \mathcal{Z}(\C)$ be the analogous braided monoidal inclusion. In both cases, we obtain for fixed $z\in\mathcal{Z}(\C)$ a map \begin{equation*}
       S_{z,\_}:\;K_0(\Rep(G))\to \K 
   \end{equation*}
   defined on the basis of $K_0(\Rep(G))$ given by simple objects
   by $t\mapsto S_{z,\tilde{\tau}(t)}$.
   The rest of the proof works as follows:\; \begin{enumerate}[label=(\roman*)]
       \item \label{ringhom1} We note that $S_{z,\_}$ is a ring homomorphism.
       \item \label{ringhom2} We note that every $g\in G$ defines a ring homomorphism $\chi_\_(g):\;K_0(\Rep(G))\to \K:\; (V, \rho_V)\mapsto \chi_{(V, \rho_V)}(g)$.
       \item \label{ringhom3} We note that all ring homomorphisms $K_0(\Rep(G))\to \K$ give a character map in the form of $\chi_\_(g)$.
   \end{enumerate}
   For \ref{ringhom1} we note that \eqref{multiplicativity} gives us that $S_{z,\_}$ is multiplicative. As we work over a fusion category, the map is additive and the additive structure is compatible with the tensor product. Furthermore, we have unitality as the double braiding of the unit object is trivial. 

   \ref{ringhom2} uses standard properties of characters.

The algebra $A=K_0(\Rep(G))\otimes_\Z \K$ is known to be
a commutative associative unital semisimple $\K$-algebra. Since
the ground field $\K$ is assumed to be algebraically closed,
it is equivalent to $\K\times\ldots \times \K$, with as many copies as
there are conjugacy classes of elements. Therefore, there are at most $|conj(G)|$-many algebra homomorphisms $A\to \K$. By \ref{ringhom2}, we have $|conj(G)|$-many ring homomorphisms $\chi_\_(g):\;K_0(\Rep(G))\to \K:\; (V, \rho_V)\mapsto \chi_{(V, \rho_V)}(g)$.
Every ring homomorphism $\chi:\; K_0(\Rep(G))\to \K$ gives a algebra homomorphism $\hat{\chi}:\;A\to \K$. The mapping $\chi\mapsto \tilde{\chi}$ is injective. Thus, every ring homomorphism is of the form in \ref{ringhom2}.

   This shows that for any $z\in \mathcal{Z}(\C)$, $S_{z,\_}$ is of the form $\chi_\_(g)$ for some $g\in G$. We conclude that we find only character values as  values $S_{z,t}$. We know that we can find, up to permutation of columns and rows, a unique invertible matrix of rank $|\Pi_0(Z_2(\C))|$ within the rows labeled by $\mathcal{T}$ in the $\Sone$-matrix of $\mathcal{Z}(\C)$. This matrix gives us the 2-categorical $\St$-matrix of $\mathcal{Z}(\Mod(\C))$. As we only have character values as matrix entries $S_{z,t}$ and as the character table of $G$ is an invertible matrix of rank $|conj(G)|=|\Pi_0(Z_2(\C))|$, we obtain the character table of $G$ as the invertible (sub)matrix we were looking for. Thus, the 2-categorical $\St$-matrix of $\mathcal{Z}(\Mod(\C))$ is the character table of $G$.
\end{proof}

\subsection{Example}\label{section32}

We now turn to a specific class of examples where the 2-categorical S-matrix can be computed directly, without recourse to Observation \ref{Smatrixcomputation1categoryviaDrinfeldcenter}. We will classify all finite semisimple indecomposable braided module categories up to equivalence and describe the Schur 
equivalence classes. 

Let $G$ be a finite abelian group and $\vect_G$ the $\Cn$-linear abelian category of finite-dimensional $G$-graded $\Cn$-vector spaces. We fix an abelian three-cocycle $(\Psi, \Omega)$ which endows $\vect_G$ with the
structure of a braided monoidal category $\C$.

To compute the 2-categorical $\St$-matrix of $\mathcal{Z}(\Mod(\C))$, we first describe the M\"uger center of $\C$.
A simple object $\so{l}\in {\vect_G^{(\Psi, \Omega)}}$ with $l\in G$ lies in the M\"uger center $\mathcal{Z}_2({\vect_G^{(\Psi, \Omega)}})$ if and only if, for all $g\in G$, $$\Omega(g,l)\Omega(l,g)=1.$$ 
The subset of such $l$ is indeed
a subgroup which we denote by $Z_2(G)$:
The unit object $\so{e}$ lies in the M\"uger center 
and the  compatibility of the braiding with the 
monoidal structure implies that $Z_2(G)$ is closed under
multiplication. Combining Example 2.1(3) and Lemma 3.2 of \cite{Ga25}, we see that the restriction $\Psi_{Z_2(G)}$ of the associator to the M\"uger center is trivial. Hence, the M\"uger center is equivalent as a symmetric monoidal category to $\vect_{Z_2(G)}^{(1,\Omega_{Z_2(G)})}$ with $\Omega_{Z_2(G)}$
the restriction of $\Omega$ to $Z_2(G)$.

Following \cite[Example 2.1]{O03} we know that the semisimple indecomposable algebras in $\vect_G^{(\Psi, \Omega)}$ are all defined via $\Psi$-trivialized
subgroups of $G$, i.e. a subgroup $H\leq G$ together 
with a 2-cochain $\mu\in C^2(H,\Cn^*)$ such that 
$d\mu=\Psi_{H}$. The cochain $\mu$ endows 
the object $A_H:=\oplus_{h\in H} \so{h}$ with the
structure of an associative separable unital algebra
$A_{(H,\mu)}$ in $\vect_G^{(\Psi, \Omega)}$.
Any finite semisimple indecomposable (right) $\vect_G^{(\Psi, \Omega)}$-module category
is equivalent to a category of left $A_{(H,\mu)}$-modules in $\vect_G^{(\Psi, \Omega)}$ for some $\Psi$-trivialized subgroup $(H,\mu)$; we denote it by
$\M_{(H,\mu)}$.

The next lemma generalizes Example \ref{regularmodulebraidings} from the classification of braidings
on the regular module category to braidings all indecomposable
module categories. It provides also explicit
examples of module categories which do not admit
the structure of a module braiding.

\begin{lemma}\label{braidingsonsubgroups}
    Let $(H,\mu)$ be a $\Psi$-trivialized subgroup of $G$, defining a module category $\M=\M_{(H,\mu)}$ over $\vect_G^{(\Psi,\Omega)}$. A module braiding on $\M$ exists if and only if $H\leq Z_2(G)$. Any module
    braiding on $\M$ is obtained by twisting the monodromy
    by an irreducible character of $G$. Every irreducible
    character can be used for such a twisting.
\end{lemma}

\begin{proof}
     Let $\so{g},\so{h}\in{\vect_G^{(\Psi, \Omega)}}$ 
     be simple objects of the monoidal category. Let $K$ be a set of representatives of the quotient $G/H$.
     Any element $k\in K$ defines a simple object $M_k$
     in the module category $\M_{(H,\mu)}$, namely the
      $A_{(H,\mu)}$-module induced from $U_k$. By definition of a module braiding $\sigma$, we have $$\sigma_{ M_k\ta \so{h},\so{g}}=\Omega(g,h)\Omega(h,g)(\sigma_{M_k,\so{g}}\ta \id_{\so{h}}).$$ 
Note that for $g\in G$ the $M_k\ta \so{g}$ is simple 
in $\M_{(H,\mu)}$; for $h\in H$, we even have 
$M_k\ta \so{h} \cong M_k$
      
       Therefore, \begin{align*}
        \langle\sigma_{M_k, \so{g}}\rangle&=\langle\sigma_{M_k\ta \so{h}, \so{g}}\rangle\\        &=\Omega(g,h)\Omega(h,g)\langle\sigma_{M_k, \so{g}}\ta\id_{\so{h}}\rangle\\
        &=\Omega(g,h)\Omega(h,g)\langle\sigma_{M_k, \so{g}}\rangle.
    \end{align*} Thus, a module braiding can only exist if $\Omega(g,h)\Omega(h,g)$ equals one for all $g\in G$,
    showing that $H$ has to be a subgroup of $Z_2(G)$. 
    
    By further inspection of the definition of a module braiding, we obtain that module braidings are given by characters of $G$: For $g\in G$ and $k\in K$ a module braiding defined by a character $\chi$ on $G$ is given by $$\sigma_{M_k, \so{g}}=\Omega(g,k)\Omega(k,g)\chi(g)\id_{M_k\otimes \so{g}}.$$ Note that this is well defined: $\Omega(g,k)\Omega(k,g)$ takes the same value on different representatives $k,k'$ of the same equivalence class in $G/H$ as $H\leq Z_2(G)$ and $\Omega(g,k)\Omega(k,g)$ is bilinear as a map $G\times G\to \Cn^*$. Examination shows that every character gives a module braiding.
\end{proof}

The last lemma together with Example \ref{regularmodulebraidings} give us all finite semisimple indecomposable braided $\vect_G^{(\Psi,\Omega)}$-module categories up to equivalence. The next step is to compute which ones of these are Schur equivalent. To this end,
we first show that Schur equivalence classes of braided module categories always have a representative with the underlying module category being the regular module category.

\begin{proposition}\label{extensionsmodulebraidingonregularmodule}
    As above, let $(H,\mu)$ be a $\Psi$-trivialized subgroup of $G$ defining a module category $\M_{(H,\mu)}$ . Assume that $\sigma$ is a module braiding on $\M_{(H,\mu)}$. Then there exists a braided module structure $\sigma'$ on the regular $\vect_G^{(\Psi, \Omega)}$-module category such that $(\M,\sigma)$ is Schur equivalent to $(\vect_G, \sigma')$.
\end{proposition}

\begin{proof}
Let $U:\;\M_{(H,\mu)}\to{\vect_G}$ be the forgetful functor. The identity
equips this functor with the structure of a module functor from
$\M_{(H,\mu)}$ to the regular $\vect_G$-module.

Let $\chi$ be the character on $G$ introduced in Lemma 
\ref{braidingsonsubgroups} that defines the module braiding on $\M_{(H,\mu)}$. 
Following Example \ref{regularmodulebraidings}, we have a module braiding
$\sigma'$ on the regular $\vect_G$-module defined by the same character:
$$\sigma'_{\so{e},\so{g}}=\chi(g)\id_{\so{e}\ta \so{g}}.$$ 
It is now straightforward to check that with this $U$ is even a braided module
functor.
\end{proof}

The previous proposition implies that we only have to consider braided structures on the regular $\vect_G^{(\Psi, \Omega)}$-module category to obtain a representative of each Schur equivalence class. The next proposition tells us which braidings on the regular $\vect_G$-module category give Schur equivalent braided module categories.

\begin{proposition}\label{modulebraidingsonregularmoduleSchurequivalence}
Let $\sigma_1,\sigma_2$ be two module braidings on the regular $\vect_G^{(\Psi, \Omega)}$-module category $\vect_G$ defined by characters $\chi_1,\chi_2$ on $G$.
\begin{enumerate}
    \item  The resulting braided module categories are Schur equivalent if and only if there exists $k\in G$ such that \begin{equation}\label{condschurequi}
        \frac{\chi_1(g)}{\chi_2(g)}=\Omega(g,k)\Omega(k,g)
    \end{equation} for all $g\in G$.
    \item The resulting braided module categories are Schur equivalent if and only if the restrictions to the M\"uger center $\mathcal{Z}_2({\vect_G^{(\Psi, \Omega)}})$ define the same braiding, i.e. for all $X\in\mathcal{Z}_2({\vect_G^{(\Psi, \Omega)}})$ in the M\"uger center of the braided monoidal
    category and all $M\in{\vect_G}$ in the module category, $$(\sigma_1)_{M,X}=(\sigma_2)_{M,X}.$$
\end{enumerate}   
\end{proposition}

\begin{proof}
\begin{enumerate}
        \item Assume there exists $k\in G$ such that $\frac{\chi_1(g)}{\chi_2(g)}=\Omega(g,k)\Omega(k,g)$. Define an endofunctor $F$ of the regular $\vect_G^{(\Psi, \Omega)}$-module category by tensoring with $\so{k}$ from the left:\; 
        $$F:\;{\vect_G}\to{\vect_G}:\;\_\mapsto \so{k}\otimes\_.$$ 
        The associator $a$ on $\vect_G^{(\Psi,\Omega)}$ endows $F$ with 
        the structure of a (left) module functor, 
        $$a_{\so{k},\so{h},\so{g}}=s_{\so{g},\so{h}}:\;F(\so{h}\ta \so{g})\to F(\so{h})\ta \so{g}.$$
    This module functor has the property of being braided, thus exhibiting
    the Schur equivalence of $({\vect_G},\sigma_1)$ and $({\vect_G},\sigma_2)$.

    For the backward direction assume that $({\vect_G},\sigma_1)$ and $({\vect_G},\sigma_2)$ are Schur equivalent braided module categories which
    is exhibited by a braided module functor 
    $$(F:\;{\vect_G}\to{\vect_G}, (s_{X,Y}:\;F(Y\ta X)\Tilde{\to} F(Y)\ta X)_{X,Y\in{\vect_G}}).$$ 
   Using direct sums in the category of (braided) module endofunctors of the regular module category, we can write
   any (braided) module functor as a direct sum of module functors mapping $\I$
   to a simple object. In the specific case of $\vect_G$, such a functor
maps simple objects to simple objects. Thus we can assume $F$ to map simple objects to simple objects.

    Hence, we find $k\in G$ such that $F(\so{e})=\so{k}$. The condition that $F$ is a braided module functor applied to $\so{g}\in \vect_G^{(\Psi, \Omega)}$ and $\so{e}$ seen as an element in the regular $\vect_G^{(\Psi, \Omega)}$-module category gives us $$\frac{\chi_1(g)}{\chi_2(g)}=\Omega(g,k)\Omega(k,g).$$
    
    \item Assume that two module braidings given by characters $\chi_1,\chi_2$ are Schur equivalent. By part (1) of this proof, we find $k\in G$ such that for all $g\in G$ $$\frac{\chi_1(g)}{\chi_2(g)}=\Omega(g,k)\Omega(k,g).$$ For $g\in Z_2(G)$, $$1=\Omega(g,k)\Omega(k,g)=\frac{\chi_1(g)}{\chi_2(g)}.$$ Thus, $\chi_1$ and $\chi_2$ coincide on all $g\in Z_2(G)$. Therefore, the restriction of the module braidings to the M\"uger center give the same braiding.

    Conversely, assume $\chi_1(g)=\chi_2(g)$ for all $g\in Z_2(G)$, i.e the module braidings are equal on $\mathcal{Z}_2(\vect_G)$. Then $$\frac{\chi_1}{\chi_2}:\;G\to \Cn^*:\;g\mapsto \frac{\chi_1(g)}{\chi_2(g)}$$ is a group homomorphism with $Z_2(G)\leq \ker(\frac{\chi_1}{\chi_2})$.
    Furthermore, $$\varphi:\;G\to \Hom(G,\Cn^*):\;g\mapsto \Omega(\_,g)\Omega(g,\_)$$ is a group homomorphism as $\Omega(\_,?)\Omega(?,\_):G\times G\to \Cn^*$ is bilinear. We know that $\ker(\varphi)=Z_2(G)$ as $\varphi(g)(h)=\Omega(h,g)\Omega(g,h)=1$ for all $h\in G$ if and only if $g\in Z_2(G)$.

    First, assume $\ker(\varphi)=Z_2(G)$ is trivial. Then $\varphi$ is injective
    and, as a map between finite sets of the same cardinality, a bijection.
 Therefore, $\frac{\chi_1}{\chi_2}\in \Ima(\varphi)$ and thus there is a $k\in G$ such that $$\Omega(\_,k)\Omega(k,\_)=\frac{\chi_1}{\chi_2}(\_).$$ By part (1) of this proposition, the two braided module categories are Schur equivalent.

    Now, let $Z_2(G)$ be non-trivial. Note that the bihomomorphism $\Omega(\_,?)\Omega(?,\_)$ factorizes through $Z_2(G)$ to a bihomomorphism $$\tilde{\beta}:\; G/Z_2(G)\times G/Z_2(G)\to \Cn^*:\; ([h],[g])\mapsto \Omega(h,g)\Omega(g,h).$$ Consider $$\varphi:\;G/Z_2(G)\to \Hom(G/Z_2(G),\Cn^*):\;[g]\mapsto \Tilde{\beta}(\_,[g]).$$ Furthermore, $\frac{\chi_1}{\chi_2}$ factors through $Z_2(G)$, too, as $Z_2(G)\leq \ker(\frac{\chi_1}{\chi_2})$. Consider, thus, $$\widetilde{\frac{\chi_1}{\chi_2}}:\;G/Z_2(G)\to \Cn^*:\;[g]\mapsto \frac{\chi_1(g)}{\chi_2(g)}.$$ Then we can proceed as in the case $Z_2(G)=\{e\}$ and obtain, by part (1) of this theorem, that the two braided module categories are Schur equivalent.
\end{enumerate}
\end{proof}

The last propositions imply the following theorem about the Schur equivalence classes of braided module categories over $\vect_G^{(\Psi, \Omega)}$.

\begin{theorem}\label{braidedmodulecatclas}
    Equivalence classes of Schur equivalent indecomposable finite semisimple braided $\vect^{(\Psi, \Omega)}_G$-module categories are in one-to-one correspondence to characters of the group $Z_2(G)$ defining the M\"uger center of $\vect^{(\Psi, \Omega)}_G$. Every Schur equivalence class has a representative with the module category being given by the regular $\vect^{(\Psi, \Omega)}_G$-module category and the braiding being defined by a character of the group $Z_2(G)$.
\end{theorem}

Now, having classified the braided module categories of $\vect_G^{(\Psi, \Omega)}$ up to Schur equivalence, we can compute the 2-categorical $\St$-matrix explicitly, thus recovering Theorem \ref{theorem:main} for pointed braided
fusion categories:

\begin{corollary}\label{StildematrixgeneralG}
    The 2-categorical $\St$-matrix of $\mathcal{Z}(\Mod(\vect_G^{(\Psi, \Omega)}))$ equals the character table of $Z_2(G)$, the group defining the M\"uger center of $\vect_G^{(\Psi, \Omega)}$.
\end{corollary}

\begin{proof}
The columns of the character table of the abelian group $Z_2(G)$  are labeled by elements of $Z_2(G)$. The same elements label isomorphism classes of simple objects of the M\"uger center of $\vect_G^{(\Psi,\Omega)}$ and 
thus columns of the $\St$-matrix. The rows of the character table are labeled by characters on $Z_2(G)$ which are,
by Theorem \ref{braidedmodulecatclas}, in one-to-one correspondence to Schur equivalence classes of $\BrMod(\vect_G^{(\Psi, \Omega)})$ and thus rows
of the $\St$-matrix. We claim that under this identification, the matrices
are identical.

For $g\in Z_2(G)$ let $\so{g}$ be the corresponding simple object in the 
M\"uger center; for a character $\chi$ on $Z_2(G)$, let $(\vect_G,\chi)$ be
the corresponding indecomposable braided module category. For simplicity,
choose $\so{e}$ as  a simple object in the regular $\vect_G^{(\Psi, \Omega)}$-module category. Following Observation \ref{halfbraidedalgebrabraidedmoduleSmatrixcomputation}, the entries of the $\St$-matrix are given by
$$\St_{(\vect_G, \chi), \so{g}}=\frac{1}{d_+(\alpha_{\so{g}})}
  \langle(\id_{\so{e}}\otimes \ev_{\so{g}})(\sigma_{\so{e},\so{g}}\otimes \alpha_{\so{g}})\circ (\id_{\so{e}}\otimes \coev_{^*\so{g}})\rangle.$$
  Using that $\vect_G^{(\Psi,\Omega)}$ has a pivotal structure inherited from
  the duality of finite-dimensional vector spaces, this simplifies to
  $$\St_{(\vect_G, \chi), {\so{g}}}= \langle\sigma_{\so{e}, {\so{g}}}\rangle= \chi(g).$$
  The 2-categorical $\St$-matrix is thus given by the character table of $Z_2(G)$.
\end{proof}

For a general braided fusion category, the situation can be more involved.
There can be Schur equivalence classes of braided module categories which 
do not have a representative with the regular module category as the 
underlying module category as the next remark show.

\begin{remark}
    Consider the case $\C=\Rep(G)$ of $\K$-linear finite-dimensional $G$-representations 
    for some finite group $G$ with the symmetric braiding inherited from vector spaces.
    By standard arguments 
    the group $\Aut_\otimes(\id_{H-mod})$ is, for any Hopf algebra $H$,  
    isomorphic to the group of central group-like elements of $H$. 
    Central group-like elements in $\Cn G$ are exactly the central elements 
    in $G$, $\Aut_\otimes(\id_{\Rep(G)})\cong Z(G)$. Thus, by Proposition \ref{braidingregular}, we have at most $|Z(G)|$-many different Schur equivalence classes of braided module categories on the regular $\C$-module category. However,
    we know that the Schur equivalence classes are in bijection to 
    simple objects of the M\"uger center and hence to simple $G$-representations
    of which there are as many as conjugacy classes in the group $G$.
    Thus, for non-abelian $G$, we do have Schur equivalence classes of braided module categories which do not have a representative of the form $({\Rep(G)}_{\Rep(G)},\sigma)$ with $\sigma$ being a module braiding on the regular module category.
\end{remark}

\begin{remark}
We present a brief remark on deformations of module braidings on
the regular module category $\C_\C$ which are, by Proposition \ref{braidingregular},
in bijection to monoidal automorphisms of the identity monoidal
functor. By \cite[Remark 4.9]{FGS24}, derivations of a monoidal functor
$F$ correspond to 1-cocycles in Davydov-Yetter cohomology $C_{DY}^1(F)$. 
For $F=\id$, they describe deformations of the identity monoidal
automorphism of $\id$. Since we work over an algebraically closed
field $\K$ of characteristic zero, \cite[Remark 4.9]{FGS24} implies that the Davydov-Yetter cohomology vanishes so that there are
no non-trivial deformations of the structure of a braided module
category on $\C_\C$.
\end{remark}
\printbibliography
\end{document}